\documentclass[11pt]{amsart}
\usepackage{amssymb}
\usepackage{amsfonts}
\usepackage{amsmath}
\usepackage{geometry}
\usepackage{graphicx}
\usepackage{amscd}

\newtheorem{theorem}{Theorem}[section]

\newtheorem{definition}[theorem]{Definition}
\newtheorem{notation}[theorem]{Notation}
\newtheorem{corollary}[theorem]{Corollary}
\newtheorem{example}[theorem]{Example}
\newtheorem{remark}[theorem]{Remark}
\numberwithin{equation}{section}

\input{tcilatex}

\begin{document}
\title{Induced Good Gradings of Structural Matrix Rings}
\author{John DeWitt}
\address{Mathematics Department, Nathan Hale High School, 11601 West Lincoln
Avenue, West Allis, WI 53227}
\email{dewittj@wawm.k12.wi.us}
\author{Kenneth L. Price}
\address{Department of Mathematics, University of Wisconsin Oshkosh, 800
Algoma Boulevard, Oshkosh, Wisconsin 54901}
\email{pricek@uwosh.edu}
\date{\today }
\maketitle

\begin{abstract}
Our approach to structural matrix rings defines them over preordered
directed graphs. A grading of a structural matrix ring is called a good
grading if its standard unit matrices are homogeneous. For a group $G$, a $G$%
-grading set is a set of arrows with the property that any assignment of
these arrows to elements of $G$ uniquely determines an induced good grading.
One of our main results is that a $G$-grading set exists for any transitive
directed graph if $G$ is a group of prime order. This extends a result of
Kelarev. However, an example of Molli Jones shows there are directed graphs
which do not have $G$-grading sets for any cyclic group $G$ of even order
greater than 2. Finally, we count the number of nonequivalent elementary
gradings by a finite group of a full matrix ring over an arbitrary field.
\end{abstract}

Many important results concerning gradings on rings and other algebraic
structures have appeared over the last two decades. Some of these are very
essential results that can serve as a guide for further study. For example,
if a grading of a Lie algebra is by a semigroup, then the semigroup must be
an abelian group (see \cite[Proposition 3.3]{Kot}). Fine gradings, good
gradings, and elementary gradings are the most fundamental types of gradings
on full matrix algebras (see \cite{BSZ} and \cite{Kot}). For finite
semigroups the good gradings of full matrix algebras can be placed in
one-to-one correspondence with simple and 0-simple primitive factors of the
semigroup (see \cite{DKvW}). For torsion-free groups, gradings with finite
support of full matrix algebras are elementary gradings (see \cite{SZ}).

More recent work has resulted in a classification of gradings for many
classical simple Lie algebras (see \cite{BKR}, \cite{EldKot}, \cite{Kot} and 
\cite{KPS}). A method to induce good gradings on Lie superalgebras, Lie
algebras, and matrix algebras using directed graphs was developed in \cite%
{DIN}, \cite{KelBook}, \cite{Kel}, and \cite{PrSz}. An analogous method to
induce good and elementary gradings on incidence rings was introduced in 
\cite{MJ} and extended in \cite{Pr}. With an increased understanding of how
to induce gradings, we can count the number of gradings of a particular
type. For full matrix algebras, fine gradings are counted in \cite{KPS} and
elementary gradings are counted in \cite{DIN}.

Throughout this note, the term \textquotedblleft ring\textquotedblright\ or
\textquotedblleft subring\textquotedblright\ will only be used for an
associative ring with unity. If we specify additive notation for a group
operation, which we will always use when we assume the group is abelian,
then the identity element will be denoted by $0$. Otherwise we use
multiplicative notation and the identity element will be denoted by $1$.

Suppose $G$ is a group and $R$ is a ring. We say $R$ is a $G$\emph{-graded
ring} if there is a direct sum $R=\bigoplus_{a\in G}R_{a}$, as a group under
the addition of $R$, such that $R_{b}R_{c}\subseteq R_{bc}$ for all $b,c\in
G $. The subgroups $R_{a}$, $a\in G$, are called the \emph{homogeneous
components} and the elements of $\cup _{a\in G}R_{a}$ are called the \emph{%
homogeneous elements}. The \emph{support} of $S$ is the set $\limfunc{Supp}%
_{G}S=\left\{ g\in G:S_{g}\neq 0\right\} $.

We define good and elementary gradings of structural matrix rings in Section %
\ref{GSMR sec}. All of our results concern a specific type of good grading,
which we call an induced good grading. As explained in Theorem \ref{Power
Principle}, for a finite abelian group $G$ the number of induced gradings
divides $\left\vert G\right\vert ^{s}$, where $s$ is the number of arrows
that are not loops of the directed graph. The construction of induced good
gradings was used in \cite{PrSz} to count induced good gradings of certain
structural matrix rings. These results also give the number of induced good
gradings of certain kinds of incidence rings (see \cite{MJ}, \cite%
{MillSpiegel}, \cite{Pr}, and \cite{PrC}).

In Section \ref{GAGG sec}, we cover an example of a directed graph $J$,
which is due to M. Jones (see \cite{MJ}), that shows a structural matrix
ring may have induced good gradings that are not elementary gradings. There
are 13 vertices in $J$. If $G$ is cyclic of order $k$, then by Theorem \ref%
{Count J} the number of induced good gradings of a structural matrix ring
over $J$ is $k^{12}$ if $k$ is odd, $2\cdot k^{12}$ if $k$ is even. In
Section \ref{GAGG sec}, we also prove Theorem \ref{Main Result}, which
extends a result of Kelarev \cite[Theorem 2]{Kel}.

Suppose $S=\bigoplus_{a\in G}S_{a}$ and $T=\bigoplus_{a\in G}T_{a}$ are $G$%
-graded rings. A \emph{homomorphism of }$G$\emph{-graded rings} is a ring
homomorphism $h:S\rightarrow T$ such that $h\left( S_{a}\right) \subseteq
T_{a}$ for all $a\in \limfunc{Supp}_{G}S$. An isomorphism that is a
homomorphism of $G$-graded rings is called an \emph{isomorphism of }$G$\emph{%
-graded rings}. In the case of matrix algebras there are gradings which are
not good gradings but are isomorphic to good gradings (see \cite[Example 1.3]%
{DIN}). Isomorphic gradings for good group gradings of incidence algebras
over partial orders have been studied by Miller and Spiegel (see \cite%
{MillSpiegel}).

Section \ref{IG sec} is dedicated to counting the number of nonequivalent
induced good gradings of a structural matrix ring. Theorem \ref{Complete
Count} provides a formula for the number of nonequivalent elementary
gradings by a finite group of a full matrix ring over a field.

\section{Graded Structural Matrix Rings\label{GSMR sec}}

For an arbitrarily chosen ring $R$, we let $M_{n}\left( R\right) $ denote
the full matrix ring on the set of square matrices over $R$. In \cite{DvW}
structural matrix rings are constructed using Boolean matrices. Instead we
may use preordered directed graphs to define structural matrix rings as in
Section 3.14 of \cite{KelBook}. However, we refer to directed edges as
\textquotedblleft arrows,\textquotedblright\ which has the advantage of not
overusing the term \textquotedblleft edges\textquotedblright\ when referring
to the underlying undirected graph (see \cite{PrSz}) or the edges of a Hasse
Diagram, as in Section \ref{GAGG sec}. It is also consistent with \cite%
{Green}, which considers gradings of path algebras that are defined using
directed graphs containing vertices and \textquotedblleft
arrows.\textquotedblright 

The directed graphs we consider have a finite number of vertices and no
repeated arrows. Loops are allowed.\ The vertex set and the arrow set of a
directed graph $D$ are denoted by $V\left( D\right) $ and $A\left( D\right) $%
, respectively. If there are $n$ vertices we may assume they are numbered so
that $V\left( D\right) =\left\{ 1,\ldots ,n\right\} $ and $A\left( D\right) $
is a subset of $V\left( D\right) ^{2}=V\left( D\right) \times V\left(
D\right) $. We often drop the parentheses and comma for any arrow $\left(
v,w\right) $ and denote it simply by $vw$. In some cases it is more
convenient to use lower case Greek letters as a notation to stand for the
arrows of a directed graph.

For $a,b\in V\left( D\right) $ with $ab\in A\left( D\right) $ we let $E_{ab}$
denote the standard unit matrix, that is, $E_{ab}$ is the $n\times n$ matrix
whose entry in row $a$ and column $b$ is 1 and all of its other entries are
0. A matrix is \textit{blocked} by $D$ if it is a linear combination of
standard matrix units which are indexed by arrows of $D$. The subset of all
blocked matrices in $M_{n}\left( R\right) $ is a free $R$-bimodule over $R$
and we denote it by $S\left( D,R\right) $.

Consider the product of two blocked standard unit matrices $B$ and $C$. We
have $BC=0$ unless $B=E_{ij}$ and $C=E_{jk}$ for some $i,j,k\in V\left(
D\right) $ such that $ij,jk\in A\left( D\right) $. Since $%
E_{ij}E_{jk}=E_{ik} $ we need $ik\in A\left( D\right) $. Therefore $S\left(
D,R\right) $ is closed under multiplication if and only if $D$ is \emph{%
transitive}. We say $D$ is \emph{reflexive} if\emph{\ }there is a loop at
every vertex of $D$. If $D$ is\emph{\ preordered}, i.e. reflexive and
transitive, then $S\left( D,R\right) $ is a ring, which is called a \emph{%
structural matrix ring}. Chapter 8 of \cite{KelBook} is devoted to gradings
of matrix rings, with gradings of structural matrix rings considered in \cite%
[Section 8.3]{KelBook}.

Let $G$ be a group. Then $S=S\left( D,R\right) $ is $G$\emph{-graded} if
there is a direct sum $S=\bigoplus_{g\in G}S_{g}$, as a group under the
addition of $S$. If $D$ is preordered, $S$ is a $G$-graded ring if $%
S_{g}S_{h}\subseteq S_{g+h}$ for all $g,h\in G$. A $G$-grading of $S\left(
D,R\right) $ is called \emph{good} if the standard unit matrices are
homogeneous. We recall some definitions and notation from \cite{PrSz}. A%
\textit{\ \emph{transitive triple}\ in }$D$ is an ordered triple of vertices
contained in 
\begin{equation*}
\limfunc{Trans}\left( D\right) =\left\{ \left( a,b,c\right) :a,b,c\in
V\left( D\right) \text{ and }ab,bc,ac\in A\left( D\right) \right\} \text{.}
\end{equation*}

A function $\Phi :A\left( D\right) \rightarrow G$ is a \emph{homomorphism}
if equation \ref{homomorphism property} holds for any $\left( a,b,c\right)
\in \limfunc{Trans}\left( D\right) $.%
\begin{equation}
\Phi \left( ab\right) \Phi \left( bc\right) =\Phi \left( ac\right)
\label{homomorphism property}
\end{equation}%
A $G$-grading of $S=S\left( D,R\right) $ is \emph{induced} by $\Phi $ if for 
$g\in G$, the homogeneous component $S_{g}$ is the span of all $E_{\alpha }$
with $\alpha \in A\left( D\right) $ such that $\Phi \left( \alpha \right) =g$%
. This is a good grading, which we call an \emph{induced good grading}. This
is not to be confused with the induced gradings defined by Bahturin and
Zaicev \cite[Definition 4.1]{BZ}, which are graded algebras formed by the
tensor product of a graded algebra with an elementary graded matrix algebra
in the case when their supports do not necessarily commute.

A homomorphism $\Phi :A\left( D\right) \rightarrow G$ is \emph{elementary}
if there are elements $g_{1},\ldots ,g_{n}\in G$ such that $\Phi \left(
ab\right) =\left( g_{a}\right) ^{-1}g_{b}$ for all $a,b\in V\left( D\right) $
such that $ab\in A\left( D\right) $. An \emph{elementary grading} of $S$ is
a grading induced by an elementary homomorphism. A subset $X$ of $A\left(
D\right) $ is a $G$\emph{-grading set for }$D$ if for every function $\phi
:X\rightarrow G$ there exists a unique homomorphism $\Phi :A\left( D\right)
\rightarrow G$ such that $\Phi |_{S}=\phi $. See \cite{PrSz} for conditions
on $D$ and $G$ that ensure all homomorphisms are elementary.

\begin{remark}
\label{elementary remark}Denote the complete graph on $n$ vertices by $%
K_{n}=\left\{ st:1\leq s,t\leq n\right\} $. For any group $G$, a $G$-grading
set for $K_{n}$ is given by $\left\{ 12,13,\ldots ,1n\right\} $ (see, for
example, \cite{DIN}, \cite{Kel} and \cite{PrSz}). Given a homomorphism $\Phi
:K_{n}\rightarrow G$, set $g_{i}=\Phi \left( 1i\right) $ for $i=1,\ldots ,n.$
For any $a,b\in G$, there is a transitive triple $\left( 1,a,b\right) $ so
equation \ref{homomorphism property} gives $\Phi \left( 1b\right) =\Phi
\left( 1a\right) \Phi \left( ab\right) $. In other words, $\Phi \left(
ab\right) =\left( g_{a}\right) ^{-1}g_{b}$. Thus the induced good gradings
of $M_{n}\left( R\right) =S\left( K_{n},R\right) $ are just the elementary
gradings.
\end{remark}

\begin{notation}
Suppose $D$ is a preordered directed graph and $G$ is a group.

\begin{enumerate}
\item $\limfunc{Hom}\left( D,G\right) $ denotes the set of homomorphisms
from $D$ to $G$.

\item $C_{G}\left( D\right) $ denotes the cardinality of the set of distinct
homomorphisms from $D$ to $G$, that is, 
\begin{equation*}
C_{G}\left( D\right) =\left\vert \limfunc{Hom}\left( D,G\right) \right\vert 
\text{.}
\end{equation*}
\end{enumerate}
\end{notation}

If all homomorphisms from $D$ to $G$ are elementary, then $C_{G}\left(
D\right) =\left\vert G\right\vert ^{n-1}$, where $n$ is the number of
vertices of $D$. Similarly, if $X$ is a $G$-grading set for $D$, then $%
C_{G}\left( D\right) =\left\vert G\right\vert ^{\left\vert X\right\vert }$.
Sufficient conditions on the directed graph for a $G$-grading set to exist
are provided in \cite{MJ}, \cite{Pr}, and \cite{PrSz}. Sufficient conditions
on the group $G$ are provided in Theorem \ref{Main Result}.

\begin{example}
\label{not an induced good grading}Suppose $G$ is a group and $R$ is a $G$%
-graded ring such that $1\in R_{1}$ and $R_{g}\neq 0$ for some $g\in
G\backslash \left\{ 1\right\} $. We may view $R$ as a good $G$-graded
structural matrix ring over $D$ with $V\left( D\right) =\left\{ 1\right\} $
and $A\left( D\right) =\left\{ \left( 1,1\right) \right\} $. There is a
homomorphism $\Phi :A\left( D\right) \rightarrow G$ given by $\Phi \left(
1,1\right) =1$. However, the grading of $R$ induced by $\Phi $ has only one
homogenous component.
\end{example}

\begin{remark}
Example \ref{not an induced good grading} shows there are good graded
structural matrix rings whose gradings are not induced by homomorphisms. We
note that this is incorrectly stated in \cite[Remark 1.2(4)]{PrSz}. However,
if $D$ is preordered and $RI\subseteq S_{1}$, where $I$ is the identity
matrix, then the good grading of $\,S=S\left( D,R\right) $ can be induced by
a homomorphism $\Phi $ defined so that $e_{ab}\in S_{\Phi \left( a,b\right)
} $ for all $a,b\in V\left( D\right) $ with $ab\in A\left( D\right) $. This
may be proved in the same way as \cite[Theorem 2]{PrC}.
\end{remark}

\begin{definition}
Let $D$ be a directed graph and let $G$ be an additively written abelian
group. For every $\alpha \in A\left( D\right) $, we let $x_{\alpha }$ denote
an unknown from $G$.

\begin{enumerate}
\item The \emph{transitive triple equations} is the homogeneous system of
linear equations of the form 
\begin{equation*}
x_{ab}-x_{ac}+x_{bc}=0
\end{equation*}%
for all $\left( a,b,c\right) \in \limfunc{Trans}\left( D\right) $ such that $%
a$, $b$, and $c$ are distinct vertices of $D$.

\item Let $r$ denote the number of transitive triple equations and let $s$
denote the number of arrows that are not loops of $D$ . The transitive
triple equations lead to the matrix equation $AX=0$ where $A$ is an $r\times
s$ integer matrix and $X$ is a vector of unknowns. We call $A$ a\textbf{\ }%
\emph{transitivity matrix} for $D$ in $G$.
\end{enumerate}
\end{definition}

\begin{theorem}
\label{Power Principle}Let $D$ be a directed graph and let $G$ be an abelian
group. If $G$ is finite, then $C_{G}\left( D\right) $ divides $\left\vert
G\right\vert ^{s}$, where $s$ is the number of arrows that are not loops of $%
D$.
\end{theorem}

\proof%
If there are no transitivity equations, then $C_{G}\left( D\right)
=\left\vert G\right\vert ^{s}$. Otherwise, $C_{G}\left( D\right) $ is equal
to the number of solutions to $AX=0$, where $A$ is a transitivity matrix for 
$D$ in $G$. The solutions to $AX=0$ form a subgroup of $G^{s}$, so $%
C_{G}\left( D\right) $ divides $\left\vert G\right\vert ^{s}$ by Lagrange's
Theorem. 
\endproof%

\section{Abelian Grading Groups \label{GAGG sec}}

\begin{figure}[th]
\begin{center}
\includegraphics[width=2.35in,height=1.53in]{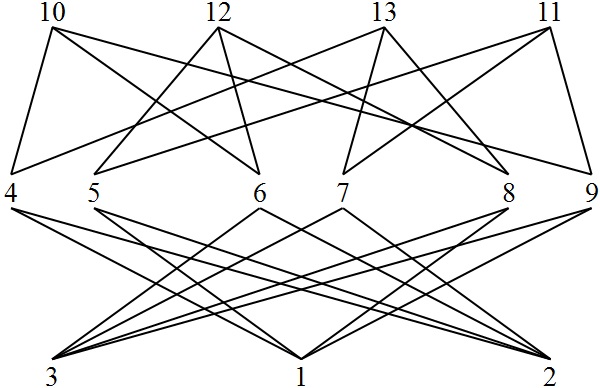}
\end{center}
\caption{Hasse Diagram for J.}
\label{Jones Digraph}
\end{figure}

Structural matrix rings can have induced good gradings that are not
elementary. We consider a particular directed graph, which we denote by $J$
since it was first studied by M. Jones (see \cite{MJ}). The Hasse diagram
for $J$ is shown in figure \ref{Jones Digraph}. The vertices are numbered
from 1 to 13. The edges all indicate arrows pointing up. Moreover, $J$ is
preordered so there are loops at every vertex in addition to all of the
arrows forced by transitivity.

\begin{notation}
For a finite group $G$, we set $\limfunc{ord}\nolimits_{k}\left( G\right)
=\left\{ g\in G:g^{k}=1\right\} $ for $k\geq 1$.
\end{notation}

\begin{theorem}
\label{Count J}Suppose $R$ is a ring and $G$ is an additively written finite
abelian group.

\begin{enumerate}
\item If $G$ has odd order, then a $G$-grading set for $J$ is given by 
\begin{equation*}
B=\left\{ \left( 3,10\right) ,\left( 3,11\right) ,\left( 3,12\right) ,\left(
3,13\right) ,\left( 2,4\right) ,\left( 2,5\right) ,\left( 3,6\right) ,\left(
2,7\right) ,\left( 3,7\right) ,\left( 3,8\right) ,\left( 1,9\right) ,\left(
3,9\right) \right\} \text{.}
\end{equation*}%
In this case, every induced good $G$-grading of $S\left( J,R\right) $ is
elementary and $C_{G}\left( J\right) =\left\vert G\right\vert ^{12}$.

\item If $G$ has order 2, then $B\cup \left\{ \left( 1,8\right) \right\} $
is a $G$-grading set for $J$ and $C_{G}\left( J\right) =2^{13}$.

\item If $G$ has even order greater than 2, then $C_{G}\left( J\right)
=\left\vert \limfunc{ord}\nolimits_{2}\left( G\right) \right\vert \cdot
\left\vert G\right\vert ^{12}$.

\item If $G$ is cyclic of even order greater than 2, then $C_{G}\left(
J\right) =2\cdot \left\vert G\right\vert ^{12}$. This is not a power of $%
\left\vert G\right\vert $, so there is no $G$-grading set for $J$.
\end{enumerate}
\end{theorem}

\begin{figure}[th]
\begin{center}
\includegraphics[width=4.17in,height=2.53in]{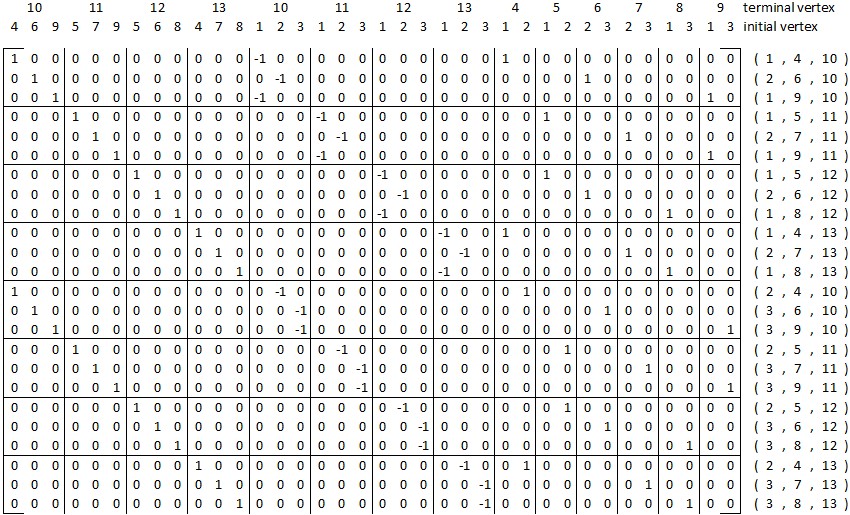}
\end{center}
\caption{Transitivity Matrix for J.}
\label{J Matrix}
\end{figure}

The transitive triple equations is a linear system of 24 equations in 36
unknowns. The transitivity matrix $A$ is shown in figure \ref{J Matrix}.
Each column of $A$ corresponds to an arrow of $J$ and every row corresponds
to a transitive triple. Every row contains two one's and a negative one,
which come from the transitivity rule. 

To solve the transitivity equations we can row-reduce $A$. Since we are
looking for solutions in an abelian group we can only use integer
row-operations. A row-echelon form of $A$ is shown in figure \ref{Reduced
Matrix}. Each row has a leading entry in a column that corresponds to a
dependent variable. The remaining variables are independent. Except for the
last row, all of the entries are either $1$, $0$, or $-1$. But in the last
row every entry is a multiple of 2. In a finite abelian group of even order
there is more than one solution to the equation $2x=0$. This leads to
another independent variable to account for these solutions.

\begin{figure}[th]
\begin{center}
\includegraphics[width=4.17in,height=2.53in]{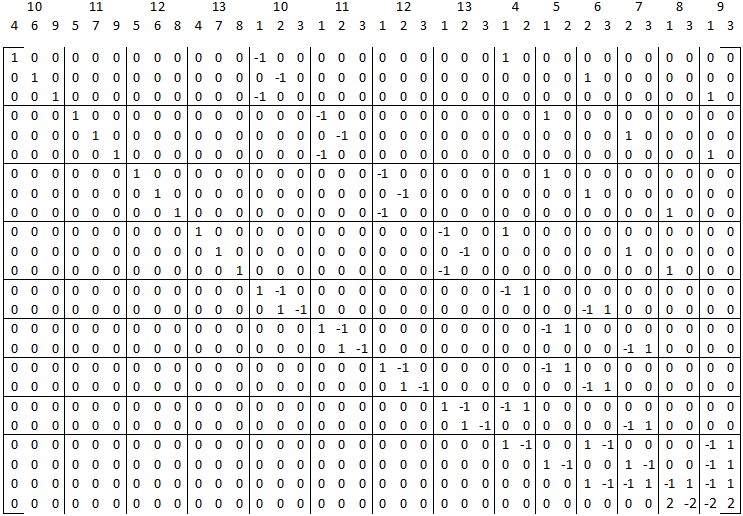}
\end{center}
\caption{Row-Reduced Matrix}
\label{Reduced Matrix}
\end{figure}

Suppose $G$ has odd order. In a parametric description the independent
variables correspond to arrows $\left( 3,10\right) $, $\left( 3,11\right) $, 
$\left( 3,12\right) $, $\left( 3,13\right) $, $\left( 2,4\right) $, $\left(
2,5\right) $, $\left( 3,6\right) $, $\left( 2,7\right) $, $\left( 3,7\right) 
$, $\left( 3,8\right) $, $\left( 1,9\right) $, and $\left( 3,9\right) $. We
find this list of arrows from the row-echelon form of the matrix. This
proves part 1, since these are the only independent variables if the order
of $G$ is odd.

If $G$ has order 2, then the last row of the row-reduced matrix yields the
equation $0=0$, so $x_{1,8}$ turns out to be an independent variable. Thus
the set $B\cup \left\{ \left( 1,8\right) \right\} $ is a $\mathbb{Z}_{2}$%
-grading set for $J$. This proves part 2.

If $G$ has even order then the last row of the row-reduced matrix yields
solutions of the form%
\begin{equation}
x_{1,8}=x_{3,8}+x_{1,9}-x_{3,9}+T  \label{Jones Eq}
\end{equation}%
where the values of the independent variable $T$ are all of the elements of $%
G$ with order 1 or 2. This proves part 3. Part 4 follows easily from part 3.%
\endproof%

\begin{remark}
Jones considered nonabelian group gradings over $J$ in \cite[Example 6]{MJ}.
Her version of equation \ref{Jones Eq} is $x_{1,8}\left( x_{3,8}\right)
^{-1}x_{1,9}\left( x_{3,9}\right) ^{-1}=g$ or $1$, where $g$ is any element
of $G$ with order 2.
\end{remark}

A result of A. V. Kelarev (see \cite[Theorem 2]{Kel}) asserts that if $%
\left\vert G\right\vert =2$, then every finite and transitive directed graph
has a $G$-grading set (which Kelarev calls a `superbasis'). We furnish a new
proof of this fact and extend the result to all finite groups of prime order.

\begin{theorem}
\label{Main Result}Suppose $D$ is a directed graph and $G$ is a finite group
of order $p$, where $p$ is prime. Then there is a $G$-grading set for $D$ of
size $s-r$, where $s$ is the number of arrows that are not loops of $D$ and $%
r$ is the rank of a transitivity matrix calculated over $\mathbb{Z}/p\mathbb{%
Z}$. Moreover, $C_{G}\left( D\right) =\left\vert G\right\vert ^{s-r}$.
\end{theorem}

\proof%
We may assume $G=\mathbb{Z}/p\mathbb{Z}$. Since $\mathbb{Z}/p\mathbb{Z}$ is
a field under addition and multiplication modulo $p$ the transitivity matrix
has a reduced row echelon form over $\mathbb{Z}/p\mathbb{Z}$. Thus there are 
$s-r$ independent variables ranging over $\mathbb{Z}/p\mathbb{Z}$ and $D$
contains a $\mathbb{Z}/p\mathbb{Z}$-grading set consisting of those arrows
that correspond to independent variables.%
\endproof%

\begin{example}
There are other directed graphs which do not contain $G$-grading sets if $G$
is cyclic of even order greater than 2, such as the one whose Hasse Diagram
is shown in Figure \ref{Second Hasse}. We do not know of any directed graphs
which do not contain a $G$-grading set when $\left\vert G\right\vert $ is odd.

\begin{figure}[th]
\begin{center}
\includegraphics[width=3.21in,height=1.71in]{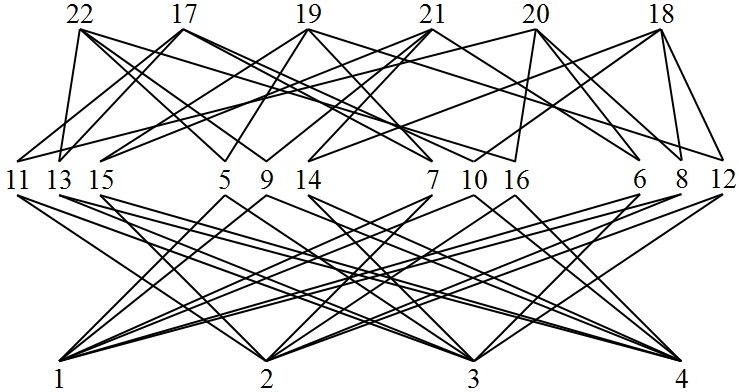}
\end{center}
\caption{Another Hasse Diagram}
\label{Second Hasse}
\end{figure}

\end{example}

\section{Equivalent Good Gradings\label{IG sec}}

\begin{definition}
Suppose $D$ is a preordered directed graph, $G$ is a group, and $R$ is a
ring. Set $S=S\left( D,R\right) $. We say that two $G$-gradings $%
S=\bigoplus_{g\in G}S_{g}$ and $S=\bigoplus_{g\in G}S_{g}^{\prime }$ of $S$
are \emph{equivalent} if there is an automorphism of $G$-graded rings $%
h:S\rightarrow S$ such that $h\left( S_{g}\right) \subseteq S_{g}^{\prime }$
for all $g\in \limfunc{Supp}_{G}S$. Let $N_{G}\left( D,R\right) $ denote the
cardinality of the set of nonequivalent induced good $G$-gradings of $%
S\left( D,R\right) $.
\end{definition}

Example \ref{DvW Example}, which is due to D\u{a}sc\u{a}lescu and van Wyk,
shows that two structural matrix rings with good gradings may be graded
isomorphic even if the underlying directed graphs are not isomorphic.

\begin{example}
\label{DvW Example}By \cite[Example 1.1]{DvW}, there is a ring $R$ such that 
\begin{equation*}
\begin{bmatrix}
R & 0 \\ 
0 & R%
\end{bmatrix}%
\cong 
\begin{bmatrix}
R & R \\ 
0 & R%
\end{bmatrix}%
\cong 
\begin{bmatrix}
R & 0 \\ 
R & R%
\end{bmatrix}%
\cong 
\begin{bmatrix}
R & R \\ 
R & R%
\end{bmatrix}%
\cong R\text{.}
\end{equation*}%
Suppose $S$ is one of the above $2\times 2$ structural matrix rings that
has\ nonzero off-diagonal entries. Then the underlying directed graph of $S$
contains an arrow $\alpha $ which is not a loop. If $g$ is a nontrivial
element of a multiplicatively written group $G$, then there is a good $G$%
-grading on $S$ such that $E_{\alpha }\in S_{g}$. Using the isomorphism $%
S\cong R$ we obtain a good grading of $R$ such that $S$ and $R$ are
isomorphic as $G$-graded rings.
\end{example}

Next we recall a few well-known terms from directed graph theory.

\begin{definition}
Let $D$, $D_{1}$, and $D_{2}$ be directed graphs.

\begin{enumerate}
\item $D_{1}$ and $D_{2}$ are \emph{isomorphic} if there is a bijection $%
\sigma :V\left( D_{1}\right) \rightarrow V\left( D_{2}\right) $ such that 
\begin{equation*}
A\left( D_{2}\right) =\left\{ \left( \sigma \left( a\right) ,\sigma \left(
b\right) \right) :a,b\in V\left( D_{1}\right) \text{ and }ab\in A\left(
D_{1}\right) \right\} \text{.}
\end{equation*}

\item The map $\sigma $ in part 1 is called an \emph{isomorphism}.

\item The \emph{arrow correspondence} is the bijection $\sigma ^{\ast
}:A\left( D_{1}\right) \rightarrow A\left( D_{2}\right) $ given by $\sigma
^{\ast }\left( ab\right) =\sigma \left( a\right) \sigma \left( b\right) $
for all $a,b\in V\left( D_{1}\right) $ such that $ab\in A\left( D_{1}\right) 
$.

\item The \emph{automorphism group} of $D$, denoted by $\limfunc{Aut}\left(
D\right) $, is the group formed by all isomorphisms from $D$ to $D$ under
composition.
\end{enumerate}
\end{definition}

\begin{remark}
Suppose $D$ is a directed graph and $G$ is a group. There is a group action
of $\limfunc{Aut}\left( D\right) $ on $\limfunc{Hom}\left( D,G\right) $
given by $\sigma .\Phi =\Phi \circ \left( \sigma ^{\ast }\right) ^{-1}$ for
any homomorphism $\Phi :D\rightarrow G$ and $\sigma \in \limfunc{Aut}\left(
D\right) $. Let $\limfunc{Hom}\left( D,G\right) _{\sigma }$ be the set of
all $\Phi \in \limfunc{Hom}\left( D,G\right) $ such that $\Phi \circ \sigma
^{\ast }=\Phi $. If $G$ is finite, then by the orbit counting formula the
number of orbits is given by%
\begin{equation}
\left\vert \limfunc{Hom}\left( D,G\right) /\limfunc{Aut}\left( D\right)
\right\vert =\frac{1}{\left\vert \limfunc{Aut}\left( D\right) \right\vert }%
\sum_{\sigma \in \limfunc{Aut}\left( D\right) }\left\vert \limfunc{Hom}%
\left( D,G\right) _{\sigma }\right\vert \text{.}  \label{Burnside's Equation}
\end{equation}
\end{remark}

The problem of counting the number of fine gradings of a full matrix algebra
over an algebraically closed field of characteristic zero was solved in \cite%
{KPS}. We will provide a formula that gives the number of nonequivalent
elementary gradings by a finite group of any full matrix ring over a field $%
\Bbbk $. By remark \ref{elementary remark}, this is equal to the number of
nonequivalent induced good $G$-gradings on the full matrix ring $S\left(
K_{n},\Bbbk \right) =M_{n}\left( \Bbbk \right) $, which is $N_{G}\left(
K_{n},\Bbbk \right) $. We list the formulas for $n\leq 6$ below. These all
follow from Theorem \ref{Complete Count}.

\begin{itemize}
\item $N_{G}\left( K_{2},\Bbbk \right) =\frac{1}{2!}\left( \left\vert
G\right\vert +\left\vert \limfunc{ord}\nolimits_{2}\left( G\right)
\right\vert \right) $

\item $N_{G}\left( K_{3},\Bbbk \right) =\frac{1}{3!}\left( \left\vert
G\right\vert ^{2}+3\left\vert G\right\vert +2\left\vert \limfunc{ord}%
\nolimits_{3}\left( G\right) \right\vert \right) $

\item $N_{G}\left( K_{4},\Bbbk \right) =\frac{1}{4!}\left( \left\vert
G\right\vert ^{3}+6\left\vert G\right\vert ^{2}+8\left\vert G\right\vert
+3\left\vert \limfunc{ord}\nolimits_{2}\left( G\right) \right\vert
\left\vert G\right\vert +6\left\vert \limfunc{ord}\nolimits_{4}\left(
G\right) \right\vert \right) $

\item $N_{G}\left( K_{5},\Bbbk \right) =\frac{1}{5!}\left( \left\vert
G\right\vert ^{4}+10\left\vert G\right\vert ^{3}+35\left\vert G\right\vert
^{2}+50\left\vert G\right\vert +24\left\vert \limfunc{ord}%
\nolimits_{5}\left( G\right) \right\vert \right) $

\item $N_{G}\left( K_{6},\Bbbk \right) =\frac{1}{6!}(\left\vert G\right\vert
^{5}+15\left\vert G\right\vert ^{4}+85\left\vert G\right\vert
^{3}+210\left\vert G\right\vert ^{2}+144\left\vert G\right\vert $\newline
$+15\left\vert \limfunc{ord}\nolimits_{2}\left( G\right) \right\vert \cdot
\left\vert G\right\vert ^{2}+90\left\vert \limfunc{ord}\nolimits_{2}\left(
G\right) \right\vert \cdot \left\vert G\right\vert +40\left\vert \limfunc{ord%
}\nolimits_{3}\left( G\right) \right\vert \cdot \left\vert G\right\vert
+120\left\vert \limfunc{ord}\nolimits_{6}\left( G\right) \right\vert )$
\end{itemize}

\begin{theorem}
\label{Elementary Lemma}Let $\Bbbk $ be a field, let $G$ be a group, and let 
$n$ be a positive integer. Suppose $S=M_{n}\left( \Bbbk \right) $ and $%
T=M_{n}\left( \Bbbk \right) $, respectively, have $G$-gradings induced by
homomorphisms $\Phi _{1}$ and $\Phi _{2}$, respectively. Then $S\ $is graded
isomorphic to $T$ if and only if there exists $\sigma \in \limfunc{Aut}%
\left( K_{n}\right) $ such that $\Phi _{2}=\Phi _{1}\circ \sigma ^{\ast }$.
\end{theorem}

\proof%
By Remark \ref{elementary remark}, $\Phi _{1}$ and $\Phi _{2}$ are
elementary homomorphisms so there exist $g_{1},\ldots ,g_{n}\in G$ and $%
h_{1},\ldots ,h_{n}\in G$ such that $\Phi _{1}\left( ab\right) =\left(
g_{a}\right) ^{-1}g_{b}$ and $\Phi _{2}\left( ab\right) =\left( h_{a}\right)
^{-1}h_{b}$ for all $a,b\in V\left( K_{n}\right) $. Suppose $S\ $is graded
isomorphic to $T$. By the same argument in the proof of \cite[Theorem 6]{BZ2}%
, there exists $g_{0}\in G$ and $\sigma \in \limfunc{Aut}\left( K_{n}\right) 
$ such that $h_{i}=g_{0}g_{\sigma \left( i\right) }$ for all $i\in V\left(
K_{n}\right) $. An easy calculation shows $\Phi _{2}\left( ab\right) =\Phi
_{1}\left( \sigma ^{\ast }\left( ab\right) \right) $ for all $a,b\in V\left(
K_{n}\right) $, as desired. The reverse implication is straightforward.%
\endproof%

We partition $S_{n}$, the permutation group on $n$ letters, following the
construction and notation of principal characteristic polynomials described
in \cite{Murnaghan}. The \emph{cycle structure} of a permutation $\sigma \in
S_{n}$ is an $n$-tuple $\alpha =\left( \alpha _{1},\ldots ,\alpha
_{n}\right) $ such that $\sigma $ is a (unique) product of disjoint $\alpha
_{1}$ cycles on one letter, $\alpha _{2}$ cycles on two letters (i.e.
transpositions), $\alpha _{3}$ cycles on three letters (i.e. ternary
cycles), and so on.

We set $\left\vert \alpha \right\vert =\alpha _{1}+\cdots +\alpha _{n}$ and $%
d\left( \alpha \right) =\gcd \left\{ i:\alpha _{i}>0\right\} $. If $\sigma
\in S_{n}$ has cycle structure $\alpha $, then $\alpha _{1}+2\alpha
_{2}+\cdots +n\alpha _{n}=n$. Thus $d\left( \alpha \right) $ divides $n$
since either $\alpha _{x}=0$ or $d\left( \alpha \right) $ divides $x$ for
every $x\leq n$. Conversely, if $d$ divides $n$, then $d\left( \alpha
\right) =d$ for $\alpha $ such that $\alpha _{d}=\frac{n}{d}$. Let $%
P_{\alpha }$ denote the number of all permutations with cycle structure $%
\alpha $. An easy derivation of the formula for $P_{\alpha }$ given in
equation \ref{coefficients}\ can be found in \cite{Murnaghan}.%
\begin{equation}
P_{\alpha }=\frac{n!}{\left( 1^{\alpha _{1}}\alpha _{1}!\right) \left(
2^{\alpha _{2}}\alpha _{2}!\right) \cdots \left( n^{\alpha _{n}}\alpha
_{n}!\right) }  \label{coefficients}
\end{equation}

\begin{theorem}
\label{Complete Count}Suppose $\Bbbk $ is a field and $G$ is a
multiplicatively written finite group. The number of nonequivalent
elementary $G$-gradings on $M_{n}\left( \Bbbk \right) $ is given by 
\begin{equation}
N_{G}\left( K_{n},\Bbbk \right) =\frac{1}{n!}\sum_{\alpha }P_{\alpha
}\left\vert \limfunc{ord}\nolimits_{d\left( \alpha \right) }\left( G\right)
\right\vert \cdot \left\vert G\right\vert ^{\left\vert \alpha \right\vert -1}%
\text{.}  \label{Matrix Formula}
\end{equation}
\end{theorem}

\proof%
Note that $V\left( K_{n}\right) =\left\{ 1,\ldots ,n\right\} $ and the
automorphism group of $K_{n}$ is $S_{n}$. By Theorem \ref{Elementary Lemma},
equivalent induced $G$-gradings on $M_{n}\left( \Bbbk \right) $ are the
orbits under the action of $S_{n}$ on $K_{n}$. Therefore, using equation \ref%
{Burnside's Equation}, we have%
\begin{equation}
N_{G}\left( K_{n},\Bbbk \right) =\frac{1}{n!}\sum_{\sigma \in
S_{n}}\left\vert \limfunc{Hom}\left( K_{n},G\right) _{\sigma }\right\vert 
\text{.}
\end{equation}

By Remark \ref{elementary remark}, $B=\left\{ 12,13,\ldots ,1n\right\} $ is
a $G$-grading set for $K_{n}$ and $C_{G}\left( K_{n}\right) =\left\vert
G\right\vert ^{n-1}$. By \cite[Theorem 2.3]{PrSz} and \cite[Proposition 4.4]%
{PrSz}, we may relabel the vertices and continue to use $B=\left\{
12,13,\ldots ,1n\right\} $ as a $G$-grading set for $K_{n}$. Moreover, two
permutations $\sigma _{1}$ and $\sigma _{2}$ have the same cycle structure $%
\alpha $ if and only if they are conjugate. In this case, $\left\vert 
\limfunc{Hom}\left( K_{n},G\right) _{\sigma _{1}}\right\vert =\left\vert 
\limfunc{Hom}\left( K_{n},G\right) _{\sigma _{2}}\right\vert $. Thus it is
enough to show equation \ref{hom formula} holds for any $\alpha $, where $%
\sigma $ is a permutation of cycle type $\alpha $. 
\begin{equation}
\left\vert \limfunc{Hom}\left( K_{n},G\right) _{\sigma }\right\vert
=\left\vert \limfunc{ord}\nolimits_{d\left( \alpha \right) }\left( G\right)
\right\vert \left\vert G\right\vert ^{\left\vert \alpha \right\vert -1}
\label{hom formula}
\end{equation}

\begin{description}
\item[Case 1] If $\alpha =\left( n,0,\ldots ,0\right) $ then $\sigma =%
\limfunc{id}$ and $\left\vert \limfunc{Hom}\left( K_{n},G\right) _{\sigma
}\right\vert =\left\vert \limfunc{Hom}\left( K_{n},G\right) \right\vert
=C_{G}\left( K_{n}\right) =\left\vert G\right\vert ^{n-1}$. Therefore,
equation \ref{hom formula} holds since $\left\vert \alpha \right\vert =n$, $%
d\left( \alpha \right) =1$ and $\limfunc{ord}\nolimits_{d\left( \alpha
\right) }\left( G\right) =\left\{ 1\right\} $.
\end{description}

In the remaining cases, $\sigma $ is a product of disjoint cycles $\sigma
_{1},\ldots ,\sigma _{k}$ of length 2 or more, with $k\geq 1$. For $1\leq
i\leq k$, let $\ell _{i}\geq 2$ denote the length of cycle $\sigma _{i}$.
Set $m_{0}=0$, $m_{1}=\ell _{1}$, $m_{2}=\ell _{1}+\ell _{2}$, $\ldots $, $%
m_{k}=\ell _{1}+\ell _{2}+\cdots +\ell _{k}$, and $m=m_{k}$. Without loss of
generality, we may assume $\sigma _{i}$ sends $m_{i-1}+1$ to $m_{i-1}+2$, $%
m_{i-1}+2$ to $m_{i-1}+3$, etc., and $m_{i}=m_{i-1}+\ell _{i}$ to $m_{i-1}+1$
for $1\leq i\leq k$. In particular, $\sigma \left( 1\right) =\sigma
_{1}\left( 1\right) =2$.

We show equation \ref{phi formula} holds for all $i,j$ such that $1<i\leq k$
and $1<j\leq \ell _{i}$.%
\begin{equation}
\Phi \left( 1,m_{i-1}+j\right) =\Phi \left( 1,2\right) ^{j-1}\Phi \left(
1,m_{i-1}+1\right)  \label{phi formula}
\end{equation}%
The transitive triple $\left( 1,2,m_{i-1}+j\right) $, the identity $\Phi
\circ \sigma ^{\ast }=\Phi $, and $\sigma ^{\ast }\left(
1,m_{i-1}+j-1\right) =\left( 2,m_{i-1}+j\right) $ are all used in the
calculation below. 
\begin{eqnarray*}
\Phi \left( 1,m_{i-1}+j\right) &=&\Phi \left( 1,2\right) \Phi \left(
2,m_{i-1}+j\right) \\
\Phi \left( 1,m_{i-1}+j\right) &=&\Phi \left( 1,2\right) \Phi \left(
1,m_{i-1}+j-1\right)
\end{eqnarray*}%
Setting $j=2$ gives $\Phi \left( 1,m_{i-1}+2\right) =\Phi \left( 1,2\right)
\Phi \left( 1,m_{i-1}+1\right) $. Setting $j=3$ and substituting gives 
\begin{equation*}
\Phi \left( 1,m_{i-1}+3\right) =\Phi \left( 1,2\right) \Phi \left(
1,m_{i-1}+2\right) =\Phi \left( 1,2\right) ^{2}\Phi \left(
1,m_{i-1}+1\right) \text{.}
\end{equation*}%
Continuing in this way, we arrive at equation \ref{phi formula}.

We have $\alpha _{2}+\cdots +\alpha _{n}=k$ since $\sigma $ is a product of
disjoint cycles $\sigma _{1},\ldots ,\sigma _{k}$ of length at least 2$.$
Moreover, $\alpha _{1}=n-m_{k}$ since $n-m_{k}$\ is the number of vertices
that are fixed by $\sigma $. This gives $\left\vert \alpha \right\vert
=\left( n-m_{k}\right) +k$. There are two remaining cases to consider.

\begin{description}
\item[Case 2] If $n=m$, then $n=m_{k}$, $\alpha _{1}=n-m_{k}=0$, $\left\vert
\alpha \right\vert =k$, and $d\left( \alpha \right) =\gcd \left\{ \ell
_{1},\ell _{2},\ldots ,\ell _{k}\right\} $. For $1<i\leq k$ and $j=\ell _{i}$%
, we show $\Phi \left( 1,2\right) \in \limfunc{ord}\nolimits_{\ell
_{i}}\left( G\right) $ using the transitive triple $\left(
1,2,m_{i-1}+1\right) $, the identity $\Phi \circ \sigma ^{\ast }=\Phi $, $%
\sigma ^{\ast }\left( 1,m_{i-1}+\ell _{i}\right) =\left( 2,m_{i-1}+1\right) $%
, and equation \ref{phi formula} in the computation below.%
\begin{eqnarray*}
\Phi \left( 1,m_{i-1}+1\right) &=&\Phi \left( 1,2\right) \Phi \left(
2,m_{i-1}+1\right) \\
\Phi \left( 1,m_{i-1}+1\right) &=&\Phi \left( 1,2\right) \Phi \left(
1,m_{i-1}+\ell _{i}\right) \\
\Phi \left( 1,m_{i-1}+1\right) &=&\Phi \left( 1,2\right) \Phi \left(
1,2\right) ^{\ell _{i}-1}\Phi \left( 1,m_{i-1}+1\right) \\
1 &=&\Phi \left( 1,2\right) ^{\ell _{i}}
\end{eqnarray*}%
Thus $\Phi \left( 1,2\right) \in \cap _{i=1}^{k}\limfunc{ord}\nolimits_{\ell
_{i}}\left( G\right) $ and $\Phi \left( 1,2\right) \in \limfunc{ord}%
\nolimits_{\gcd \left\{ \ell _{1},\ell _{2},\ldots ,\ell _{k}\right\}
}\left( G\right) $. By equation \ref{phi formula}, $\Phi $ is completely
determined by $\Phi \left( 1,2\right) \in \limfunc{ord}\nolimits_{d\left(
\alpha \right) }\left( G\right) $ and $\Phi \left( 1,m_{1}+1\right) ,\ldots
,\Phi \left( 1,m_{k-1}+1\right) \in G$. Equation \ref{hom formula} holds in
this case with 
\begin{equation*}
\left\vert \limfunc{Hom}\left( K_{n},G\right) _{\sigma }\right\vert
=\left\vert \limfunc{ord}\nolimits_{d\left( \alpha \right) }\left( G\right)
\right\vert \left\vert G\right\vert ^{k-1}=\left\vert \limfunc{ord}%
\nolimits_{d\left( \alpha \right) }\left( G\right) \right\vert \left\vert
G\right\vert ^{\left\vert \alpha \right\vert -1}\text{.}
\end{equation*}

\item[Case 3] If $n>m$, then $n>m_{k}$, $\alpha _{1}\neq 0$, $d\left( \alpha
\right) =1$, and $\limfunc{ord}\nolimits_{d\left( \alpha \right) }\left(
G\right) =\left\{ 1\right\} $. We use the transitive triple $\left(
1,2,m_{k}+1\right) $, the identity $\Phi \circ \sigma ^{\ast }=\Phi $, and $%
\sigma ^{\ast }\left( 1,m_{k}+1\right) =\left( 2,m_{k}+1\right) $ in the
calculation below. 
\begin{eqnarray*}
\Phi \left( 1,m_{k}+1\right) &=&\Phi \left( 1,2\right) \Phi \left(
2,m_{k}+1\right) \\
\Phi \left( 1,m_{k}+1\right) &=&\Phi \left( 1,2\right) \Phi \left(
1,m_{k}+1\right) \\
1 &=&\Phi \left( 1,2\right)
\end{eqnarray*}%
By equation \ref{phi formula}, $\Phi $ is completely determined by $\Phi
\left( 1,m_{1}+1\right) ,\ldots ,\Phi \left( 1,m_{k-1}+1\right) \in G$ and $%
\Phi \left( 1,m_{k}+1\right) ,$ $\Phi \left( 1,m_{k}+2\right) $, $\ldots $, $%
\Phi \left( 1,n\right) \in G$. Equation \ref{hom formula} holds in this case
with 
\begin{equation*}
\left\vert \limfunc{Hom}\left( K_{n},G\right) _{\sigma }\right\vert
=\left\vert G\right\vert ^{k-1}\left\vert G\right\vert ^{n-m_{k}}=\left\vert 
\limfunc{ord}\nolimits_{d\left( \alpha \right) }\left( G\right) \right\vert
\left\vert G\right\vert ^{\left\vert \alpha \right\vert -1}\text{.}
\end{equation*}
\end{description}

\endproof%

\begin{remark}
Suppose $n$ and $\left\vert G\right\vert $ are relatively prime. Recall $%
d\left( \alpha \right) $ divides $n$ for every $\alpha $. Thus $d\left(
\alpha \right) $ and $\left\vert G\right\vert $ are relatively prime, which
means $\limfunc{ord}\nolimits_{d\left( \alpha \right) }\left( G\right)
=\left\{ 1\right\} $. In this case formula \ref{Matrix Formula} simplifies
to formula \ref{Matrix Formula 2}. 
\begin{equation}
N_{G}\left( K_{n},\Bbbk \right) =\frac{1}{n!}\sum_{\alpha }P_{\alpha
}\left\vert G\right\vert ^{\left\vert \alpha \right\vert -1}\text{.}
\label{Matrix Formula 2}
\end{equation}%
Formula \ref{principal characteristic} shows Murnaghan's notation for $q_{n}$%
, the principal characteristic of $S_{n}$, where $s_{1},s_{2},\ldots ,s_{n}$
are indeterminates (see \cite{Murnaghan}).%
\begin{equation}
q_{n}\left( s_{1},s_{2},\ldots ,s_{n}\right) =\sum_{\alpha }\frac{1}{\alpha
_{1}!\alpha _{2}!\cdots \alpha _{n}!}\left( \frac{s_{1}}{1}\right) ^{\alpha
_{1}}\left( \frac{s_{2}}{2}\right) ^{\alpha _{2}}\cdots \left( \frac{s_{n}}{n%
}\right) ^{\alpha _{n}}  \label{principal characteristic}
\end{equation}%
Then formula \ref{Matrix Formula 2} becomes $N_{G}\left( K_{n},\Bbbk \right)
=q_{n}\left( \left\vert G\right\vert ,\ldots ,\left\vert G\right\vert
\right) \cdot \left\vert G\right\vert ^{-1}$.
\end{remark}

As an application of Theorem \ref{Complete Count}, we offer Corollary \ref%
{Finite Groups corollary}, which is a well-known result for finite groups.

\begin{corollary}
\label{Finite Groups corollary}Suppose $G$ is a finite group. Then%
\begin{equation*}
\left\vert \limfunc{ord}\nolimits_{p}\left( G\right) \right\vert \equiv
\left\vert G\right\vert ^{p-1}\text{ }\left( \func{mod}p\right)
\end{equation*}%
for any prime number $p$.
\end{corollary}

\proof%
Setting $n=p$ and multiplying by $p!$ on both sides of equation \ref{Matrix
Formula} gives 
\begin{equation*}
\left( p!\right) N_{G}\left( K_{p},\Bbbk \right) =\sum_{\alpha }P_{\alpha
}\left\vert \limfunc{ord}\nolimits_{d\left( \alpha \right) }\left( G\right)
\right\vert \cdot \left\vert G\right\vert ^{\left\vert \alpha \right\vert -1}%
\text{.}
\end{equation*}%
We show that reducing mod $p$ yields the result since all of the terms
disappear except in cases 1 and 2 below.

\begin{description}
\item[Case 1] If $\alpha =\left( p,0,\ldots ,0\right) $, then $P_{\alpha }=1$
by formula \ref{coefficients} and $\left\vert \limfunc{ord}%
\nolimits_{d\left( \alpha \right) }\left( G\right) \right\vert \cdot
\left\vert G\right\vert ^{\left\vert \alpha \right\vert -1}=1\cdot
\left\vert G\right\vert ^{p-1}$.

\item[Case 2] If $\alpha =\left( 0,\ldots ,0,1\right) $, then, by formula %
\ref{coefficients}, $P_{\alpha }=\left( p-1\right) !$ and $\left( p-1\right)
!\equiv -1$ $\left( \func{mod}p\right) $ by Wilson's Theorem. Moreover, $%
\left\vert \limfunc{ord}\nolimits_{d\left( \alpha \right) }\left( G\right)
\right\vert \cdot \left\vert G\right\vert ^{\left\vert \alpha \right\vert
-1}=\left\vert \limfunc{ord}\nolimits_{p}\left( G\right) \right\vert \cdot 1$%
.

\item[Case 3] For $\alpha \notin \left\{ \left( p,0,\ldots ,0\right) ,\left(
0,\ldots ,0,1\right) \right\} $, we have $\alpha _{1},\ldots ,\alpha
_{p-1}<p $ and $\alpha _{p}=0$ so none of the terms in the denominator of
formula \ref{coefficients} is a multiple of $p$. Thus $P_{\alpha }\equiv 0$ $%
\left( \func{mod}p\right) $.
\end{description}

\endproof%

\paragraph{Acknowledgements.}

The authors express their gratitude to the referee. One of the authors,
Price, dedicates his contributions to Dr. Jeffrey Bergen, who recently
passed away.

\end{document}